\documentclass[11pt]{article}
\usepackage{amsthm, amsfonts, url}
\usepackage{fullpage}
\usepackage{amscd,amsmath,amssymb}

\newtheorem{definition}{Definition}
\newtheorem{theorem}{Theorem}
\newtheorem{lemma}{Lemma}
\newtheorem{proposition}{Proposition}

\def\al{\alpha}

\def\La{\Lambda}
\def\la{\lambda}

\def\kappa{\varkappa}

\def\C{{\mathbb C}}
\def\y{{\mathbb Y}}
\def\Z{{\mathbb Z}}
\def\yn{{\mathbb Y}_n}
\def\sn{{\mathfrak S}_n}
\def\si{{\mathfrak S}_{\mathbb N}}
\def\gn{{\mathbb C}[\sn]}

\def\gzn{\operatorname{GZ}_n}

\def\Comm{\operatorname{Comm}}

\def\N{{\mathbb N}}
\def\k{{\mathbf k}}

\def\S{{\frak S}}

\def\Prob{\operatorname{Prob}}
\def\Ind{\operatorname{Ind}}

\def\tr{\operatorname{tr}}

\def\be{\begin{equation}}
\def\ee{\end{equation}}

\urldef\vershik\url{vershik@pdmi.ras.ru}
\urldef\natalia\url{natalia@pdmi.ras.ru}

\author {N.~V.~Tsilevich\thanks{%
St.~Petersburg Department of Steklov Institute of Mathematics.
E-mail: \natalia, \vershik. Supported by the grants
CRDF RUM1-2622-ST-04, INTAS 03-51-5018,
RFBR 05-01-00899, and NSh-4329.2006.1.%
} \and A.~M.~Vershik\footnotemark[1]}
\title {Induced representations of the infinite symmetric group}

\date{May 3, 2006}

\begin{document}
\maketitle

\begin{abstract}
We study the representations of the infinite symmetric group induced
from the identity representations of Young subgroups. It turns out
that such induced representations
can be either of type~I or of type~II. Each Young subgroup corresponds
to a partition of the set of positive integers;
depending on the sizes of blocks of this partition, we divide Young subgroups
into two classes:
large and small subgroups.
The first class gives representations of type I, in particular, irreducible
representations. The most part of Young subgroups of
the second class give representations of type~II and, in particular,
von Neumann factors of
type II. We present a number of various examples. The main
problem is to find the so-called spectral measure of the induced representation.
The complete solution of this problem is given for two-block Young subgroups
and subgroups with infinitely many singletons and finitely many finite blocks
of length greater than one.
\end{abstract}

\tableofcontents

\section{Introduction}

In the classical representation theory of
the symmetric groups, representations induced from Young subgroups
(i.e., subgroups that leave some partition fixed)
play a very prominent role (see, e.g., \cite{Fulton}).
The irreducible components of such representations contain
all irreducible representations, and in the classical approach,
this allows one to establish a connection between Young diagrams
and irreducible representations. Though at present there is
an alternative approach to establishing this correspondence
(avoiding induction from Young subgroups), nevertheless, the traditional
problems  concerning induction for representations of the infinite symmetric
group are of independent importance. In this paper, we start the
systematic study of representations of
the infinite symmetric group $\si$ induced from infinite Young subgroups.
Let us briefly present the main results.

One can roughly divide partitions and the corresponding Young
subgroups into two classes: {\it large} partitions, which have
finitely many finite blocks and an arbitrary number of infinite blocks
(in this case, there is at least one infinite block), and {\it small}
partitions, which have infinitely many finite blocks (without any
assumption on the number of infinite blocks).
In contrast to the case of finite symmetric groups, induction from
the identity representations of a  Young subgroup to $\si$ often
gives an irreducible representation --- this is the case for
large Young subgroups corresponding to partitions with at most
one finite block. This can be proved using the well-known
Mackey criterion of irreducibility, or by direct arguments.
Induction from an arbitrary large subgroup gives a type I
representation, which, as is well known, has a unique decomposition
into irreducible ones, and we describe this decomposition explicitly. It
turns out to be finite, and, in general, its irreducible components are
not induced representations.

In contrast to this case, induction from small Young subgroups
gives type II representations. There is a well-developed theory
of representations of $\si$ of type II$_1$ with finite trace, but
here we have another situation, which was not paid much attention to.
Namely, the factor generated by the representation operators
is of type II$_\infty$, so that it has no finite trace, but
its commutant is of type II$_1$. Note that even in the case
where the factor generated by the representation operators is of type II$_1$, it may happen
that it has no finite trace, because the cyclic vector of
the representation is not cyclic for the commutant.
In \cite{strange}, such representations were called {\it strange}.
The additional invariant of such representations
is the so-called {\it coupling constant}.

It seems that the classification of such representations (up to
equivalence) was not considered; it is not even known whether it is
simpler than the classification of irreducible representations,
which, as is well known, is a wild problem for groups that are not
virtually commutative, such as $\si$.
A very interesting problem is whether each irreducible
representation of $\si$ appears in the decompositions of
induced representations (more exactly, is weakly contained
in induced representations), as is the case for the finite
symmetric groups. It should be noted that there are few papers devoted to induced
representations of the infinite symmetric group
(see, e.g., papers by Binder \cite{Binder1, Binder2, Binder3}, Obata
\cite{Obata1, Obata2}, and Hirai
\cite{Hirai1, Hirai2}). In particular, \cite{Binder1}
contains the irreducibility criterion for induction from Young
subgroups (our Theorem~\ref{th:typeI}(b)).

Models of representations of the infinite symmetric group $\si$ can be roughly
divided into two classes (see \cite{AAM}): {\it substitutional} models
and {\it spectral} models. The first class contains models in which
a group $G$ (not necessarily $\si$) acts by substitutions on
some $G$-space. Induced (quasi-regular) representations belong to
this type. The notion of a spectral model uses specific properties of
the infinite symmetric group.
The spectral analysis of a representation of a finite
symmetric group means its decomposition into
irreducible components indexed
by Young diagrams (with some multiplicities). In a more
rigorous style, this means that we diagonalize the image of the
group algebra in the representation with respect to the so-called
commutative Gelfand--Tsetlin algebra (see \cite{VO}). The
group algebra of $\si$ also has a
a distinguished maximal commutative subalgebra, the so-called
{\it Gelfand--Tsetlin algebra}, whose spectrum
is the space of infinite Young tableaux. The space of every cyclic unitary representation
of $\si$ can be {\it diagonalized} with respect to this  algebra.
This gives a spectral model of this representation, which is realized
in the Hilbert space $L^2_\mu(T,H)$, where $T$ is the space of infinite
Young tableaux, $\mu$ is a measure on $T$ (the spectral measure
of a cyclic vector), and $H$ is an auxiliary Hilbert space. Thus two problems
arise, which are similar to problems of Fourier theory:

1) (Direct problem). To find a spectral model of a representation
given in a substitutional realization.

2) (Inverse problem). To find a substitutional realization of
a representation given in a spectral form.

We consider these two problems as the main problems of this theory, and
its main analytic and
probabilistic component is
the analysis of spectral measures; see, e.g.,
the authors' papers \cite{AAM, VTs, 325}.
In this paper, we consider the first
problem for some classes of induced representations.
The first result in this direction, obtained in \cite{VTs} and specified below,
is that the representations induced from Young
subgroups corresponding to partitions into two infinite blocks
are simple and irreducible, and their spectral measures are Markov
measures. Another nontrivial example considered in this paper
is the class of induced representations corresponding to partitions
with infinitely many singletons and finitely many finite blocks
of length greater than one. In this case,
the spectral measure is a convex combination of conditional Plancherel
measures. A more precise spectral analysis of induced
representations of $\si$ will be presented elsewhere.

The paper is organized as follows. In Sec.~\ref{sec:not}, we
introduce necessary definitions and notation
related to induced representations of $\si$
we are going to consider. Section~\ref{sec:typeI}
is devoted to induced representations of type~I; namely, we give a condition
under which the induced representation is of type~I, describe the decomposition
of such a representation into irreducible components, give an
irreducibility criterion, and consider several classes of examples.
In Sec.~\ref{sec:factor}, in a similar way we deal with
induced representations of type~II: we give a condition
under which the induced representation is of type~II,
describe the central decomposition of such a representation into factors, give a criterion
of its being a factor, and consider several examples.
Finally, in Sec.~\ref{sec:spectral}, we present two
examples of the spectral analysis of the induced representations
corresponding to  partitions with two infinite blocks and to
partitions with infinitely many singletons and finitely many finite blocks
of length greater than one.

\section{Young subgroups and induced representations}
\label{sec:not}

We denote by $\sn$ the symmetric group of degree $n$.

The irreducible representations of the symmetric group
$\sn$ are indexed by the set
$\yn$ of Young diagrams with $n$ cells. Let
$\pi_\la$ be the irreducible unitary representation of
$\sn$ corresponding to a diagram
$\la\in\yn$, and let $\dim\la$ be
the dimension of $\pi_\la$.
The branching of irreducible representations of the symmetric groups
is described by the Young graph $\y$. The set of vertices of the $\Z_+$-graded
graph $\y$ is
$\cup_n\yn$, and two vertices $\mu\in{\mathbb Y}_{n-1}$ and $\la\in\yn$
are joined by an edge if and only if $\mu\subset\la$.
By definition, the zero level ${\mathbb Y}_{0}$ consists of
the empty diagram $\emptyset$.

Now let
$\si=\cup_{n=1}^\infty\sn=\varinjlim\sn$ be the infinite symmetric group
with the fixed structure of an inductive limit of finite groups.

 Consider an arbitrary partition $\Pi=(A_1,A_2,{\ldots})$ of the
set of positive integers $\N$ into disjoint subsets $A_1,A_2,{\ldots}$.
The corresponding
{\it Young subgroup} of the infinite symmetric group
$\si$ is $\S_\Pi={\frak S}_{A_1}\times {\frak S}_{A_2}\times{\ldots}$,
where ${\frak S}_A$ is the group of all finite permutations
of the elements of a set $A$.

\begin{definition}
Let $\Pi=(A_1,A_2,{\ldots})$ be a partition of $\N$.
For $i=1,2,{\ldots} ,\infty$,
denote by $k_i\in\N\cup\{\infty\}$ the number of $j$ such that
$|A_j|=i$.
We say that the partition $\Pi$ and the corresponding
Young subgroup $\S_\Pi$ are {\em of type}
$\k=(\infty^{k_\infty},1^{k_1},2^{k_2},{\ldots})$.
\end{definition}

It turns out that it is natural to divide partitions and the corresponding
Young subgroups into two classes. Namely, we introduce the following
definition.

\begin{definition}
A partition $\Pi$ of $\N$ and the corresponding Young subgroup $\S_\Pi$
is called {\em large} if it has
finitely many finite blocks
(and, necessarily, at least one infinite block). Otherwise, i.e.,
if $\Pi$ has infinitely many finite blocks, it is called  {\em small}.
\end{definition}

Our purpose is to investigate the representations
$I_\Pi=\Ind_{\S_\Pi}^{\si}1_{\S_\Pi}$
of the infinite symmetric
group $\si$ induced from the identity representations of Young
subgroups.

Note that two Young subgroups $\S_{\Pi_1}$ and $\S_{\Pi_2}$
of the same type can be sent to each other by an automorphism
of $\si$.
If the partitions $\Pi_1$ and $\Pi_2$
are ``tail-equivalent'' (i.e., can be obtained from each other
by a finite permutation of $\N$)
then this automorphism is an inner automorphism of $\si$, so that
the subgroups $\S_{\Pi_1}$ and $\S_{\Pi_2}$ are conjugate in $\si$, and
the representations $I_{\Pi_1}$ and $I_{\Pi_2}$ are obviously equivalent.
In fact, the following lemma holds.

\begin{lemma}
The representations $I_{\Pi_1}$ and $I_{\Pi_2}$ are equivalent
if and only if the Young subgroups $\S_{\Pi_1}$ and $\S_{\Pi_2}$
are conjugate in $\si$ (that is, the partitions $\Pi_1$ and $\Pi_2$ are
tail-equivalent).
\end{lemma}
\begin{proof}
The ``if'' part is obvious. Let us prove the ``only if'' part. Assume that
$I_{\Pi_1}$ and $I_{\Pi_2}$ are equivalent, and let $A$ be the
corresponding intertwining operator.
As in the proof of Theorem~\ref{th:typeI}, one can show
that this operator is determined by a function $\al$ defined on the left
cosets $\si/\S_{\Pi_2}$ and constant on the orbit of the left action
of $\S_{\Pi_1}$, that is, on left cosets lying in the same double coset
$\S_{\Pi_1}\backslash\si/\S_{\Pi_2}$, so that $\al$ is concentrated
on finite orbits, i.e., on
double cosets that decompose into the union of finitely many left cosets.
Then it is not difficult to show that there are no such finite orbits unless
the partitions $\Pi_1$ and $\Pi_2$ differ by at most finitely many blocks,
so that the lemma follows from the fact that for a finite symmetric group,
the representations induced from the identity representations of
two Young subgroups are equivalent if and only if these subgroups
are conjugate.
\end{proof}

Note that for irreducible induced representations, this assertion
follows from a result of Mackey
(see, e.g., \cite[Corollary~3, p.~158]{Mackey}).

We see that partitions of the same type $\k$ can lead to a continuum of
nonequivalent induced representations.
However, as we will see below,
the rough properties of these representations (such as being
irreducible or being a factor) depend only on the type $\k$.

Note that, in general, the properties of induced representations
$\Ind_H^G{1_H}$ are not continuous with respect to $H$.
For example, consider the simplest case --- an increasing sequence
of partitions of type $(\infty,n)$ (one infinite block
and one finite block of size $n$). For a given $n$, all such Young
subgroups are conjugate, so that the corresponding induced representations
are equivalent. Hence the limits of sequences of invariants of induced
representations also coincide. But such a sequence of partitions
can converge to different types of partitions: in the limit, we have either the trivial partition into
one block and hence the identity representation, or
a partition with two infinite blocks, in which case the possible
limits are a series of two-block representations, including a continuum
of nonequivalent ones. At the same time, if we
consider an increasing sequence of partitions with {\it fixed cyclic vectors}
and embed them into one another preserving the cyclic vectors, then
the matrix elements will be obviously continuous.

Thus, in what follows, we denote by $\xi$ the (normalized)
{\it distinguished cyclic vector}
of the representation $I_\Pi$.
If we consider the realization of $I_\Pi$ in the
$l^2$ space over the homogeneous space $\si/\S_\Pi$, then
$\xi$ is the delta function at the coset $\S_\Pi\in\si/\S_\Pi$.
Let $\Pi_n=\Pi\cap\{1,{\ldots} ,n\}$, denote by
$I_{\Pi_n}$ the representation of the finite symmetric group $\sn$
induced from the identity representation of
the Young subgroup $\S_{\Pi_n}\subset\sn$, and let
$\xi_n$ be the distinguished cyclic vector in $I_{\Pi_n}$.
Then {\it $I_\Pi$ is the inductive limit of $I_{\Pi_n}$, and the corresponding
embeddings preserve the distinguished cyclic vectors}.
Note also that $I_{\Pi_n}$ can be identified with
the restriction of $I_\Pi$ to  $\sn$
in the cyclic hull $\sn\xi$ of $\xi$; then $\xi_n$
is identified with $\xi$.
\section{Induced representations of type I}
\label{sec:typeI}

\subsection{Large Young subgroups lead to type I representations}

\begin{theorem}
\label{th:typeI}
{\rm (a)} The representation of the infinite symmetric group $\si$ induced
from the identity representation of a large Young subgroup of type
$\k=(\infty^{k_\infty},1^{k_1},2^{k_2},{\ldots} )$ with finitely many finite blocks
(that is, with $k_1+k_2+{\ldots}<\infty$) is of type I and decomposes into
a finite sum of irreducible representations.

{\rm (b)} The induced representation $I_\Pi$
of type $(\infty^{k_\infty},1^{k_1},2^{k_2},{\ldots})$
is irreducible if and only if the partition $\Pi$
contains at most one finite element, that is, $k_1+k_2+{\ldots} \le 1$.
\end{theorem}

\begin{proof}
(a) Let $\Pi=(A_1,A_2,{\ldots} )$ be a partition of type $\k$ with
finitely many finite blocks. For convenience, assume that
the blocks $A_1,{\ldots} ,A_n$ of $\Pi$ are finite
and the remaining blocks are infinite.
For brevity, denote by $H=\S_\Pi$ the corresponding Young subgroup, and
let $\S_{\rm fin}=
\S_{\rm fin}(\Pi)={\frak S}_{A_1}\times{\ldots} \times{\frak S}_{A_n}$
be the {\it finite} Young subgroup corresponding to the finite blocks of $\Pi$.

By definition, the induced representation $I_\Pi$ acts in the
space $l^2(X)$ of square-summable functions on the (countable) homogeneous
space $X=\si/H$.
Denote by $R(I_\Pi)$ the space of intertwining operators for the
representation $I_\Pi$, and
let $A\in R(I_\Pi)$. The basis of $l^2(X)$ consists of the
delta functions at the left cosets $sH\in X$ of $H$ in $\si$; by abuse of notation,
we will denote a basis element by the same symbol as the corresponding coset.
Then the operator $A$ is determined
by a matrix $(a_{sH,tH})_{sH,tH\in X}$ in this basis.
Let
$F\subset H\backslash G/H$ be the set
consisting of those double cosets $HgH$ that are unions
of finitely many left cosets of $H$; and let
$\Comm_{\si}(H)$ be the group, called the commensurator of $H$ in $\si$,
consisting of $g\in\si$ such that $HgH\in F$.
Denote by ${\mathbf1}_x$ the
characteristic function of a double coset $x\in F$.
Using the standard arguments, which go back to Mackey
and his successors, one can show that
the set of intertwining operators
$R(I_\Pi)$ is generated by the operators $A_x$, $x\in F$,
with matrix elements
$a_{sH,tH}={\mathbf1}_x(s^{-1}tH)$.

In our case, it is easy to see that $g\in\Comm_H^{\si}$ if and only if
$g$ leaves invariant all infinite blocks of the partition $\Pi$, i.e.,
$g\in \S_{\rm fin}$.
It follows that $R(I_\Pi)$ is isomorphic to the space $R(I_{\rm fin})$
of intertwining operators for the representation $I_{\rm fin}$
of the finite symmetric group ${\frak S}_N$, where $N=k_1+k_2+{\ldots}$,
induced from the identity representation of the Young subgroup
$\S_{\rm fin}$. Denote by
$P_i\in R(I_{\rm fin})$ the projections to the irreducible components
of $I_{\rm fin}$. Then their images $\tilde P_i$ in $R(I_\Pi)$
determine the finite decomposition of $I_\Pi$ into irreducible components,
and assertion (a) follows.

(b) It is easy to see that the condition $k_1+k_2+{\ldots} \le 1$
is equivalent to $\Comm_{\si}(H)=H$, and the latter condition means that
there are no nonscalar intertwining operators.
\end{proof}

\noindent{\bf Remark.}
The claim (b) of this theorem was proved in \cite{Binder1}.
It also easily follows
from the irreducibility criterion for induced
representations of discrete groups due to Mackey (see, e.g.,
\cite[Corollary~2, p.~158]{Mackey}).

\subsection{Decomposition of type I representations into irreducible components}
Let $\S_\la={\frak S}_{A_1}\times{\ldots} \times{\frak S}_{A_m}$
be a Young subgroup of a finite symmetric group
${\frak S}_N$ associated with a Young diagram $\la$
(which means that the lengths $|A_k|$ of blocks form the diagram $\la$). Then
the decomposition of
the representation of ${\frak S}_N$ induced from the Young subgroup $\S_\la$
into irreducible components is given by the following formula
(see \cite{Fulton}
and also \cite{V}):
\begin{equation}
\label{ind}
\Ind_{\S_\la}^{{\frak S}_N} 1_{\S_\la}=\sum_{\mu\unrhd\la}K_{\mu,\la}\pi_\mu,
\end{equation}
where $\unrhd$ is the dominance ordering on partitions
(see \cite[Sec.~I.1]{Mac}) and $K_{\mu,\la}$ are the Kostka numbers
(see \cite[Sec.~I.6]{Mac}). (Recall that $\mu\unrhd\la$ if and only if
$\mu_1+\mu_2+{\ldots} +\mu_i\ge\la_1+\la_2+{\ldots} +\la_i$ for all $i\ge1$,
and $K_{\mu,\la}$ is the number of semistandard Young tableaux of shape $\mu$
and weight $\la$.)

Combining these observations with the proof of Theorem~\ref{th:typeI},
we obtain the following theorem.

\begin{theorem}
Let $\Pi$ be a partition of $\N$ with finitely many finite blocks.
Denote by $\la=\la(\Pi)=1^{k_1}2^{k_2}{\ldots}$ the Young diagram
of size $N=k_1+k_2+{\ldots}$
formed by the lengths of the finite blocks of $\Pi$. Then
\begin{equation}
\label{dec}
I_\Pi=\sum_{|\mu|=N,\,\mu\unrhd\la}K_{\mu,\la}\pi_\mu^{k_\infty},
\end{equation}
where $\pi_\mu^{k_\infty}$ is the irreducible representation of
$\si$ corresponding
to the irreducible component $\pi_\mu$ of $I_{\rm fin}$ as described in the
last paragraph of the proof of Theorem{\rm~\ref{th:typeI}(a)}.
\end{theorem}

The representation  $\pi_\mu^{k_\infty}$ can be explicitly described as follows.
Let $m_i=|A_i|$.
It is not difficult to see that the homogeneous space $X_\Pi=\si/\S_\Pi$
can be realized as the space of {\it ordered} partitions of $\N$
into disjoint sets of given cardinalities
\begin{equation}
\label{hom}
r=(R_1,R_2,{\ldots} ),\qquad |R_i|=m_i,
\end{equation}
tail-equivalent to the original partition $\Pi$
with the natural action of
$\si$  by permutations.
Moreover,  the subset of all such partitions with
fixed infinite blocks $R_{n+1}, R_{n+2},{\ldots} $
can be identified with the homogeneous space ${\frak S}_N/\S_{\rm fin}$,
and the corresponding set of functions can be identified with the space
of the representation
$I_{\rm fin}=\Ind_{\S_{\rm fin}}^{{\frak S}_N} 1$.
Let $P_\mu$ be the projection to an irreducible component $\pi_\mu$ of $I_{\rm fin}$.
Then the representation $\pi_\mu^{k_\infty}$ acts in the subspace obtained
by applying the projection $P_\mu$ to the first ``coordinates''
$(R_1,{\ldots} ,R_n)$ and leaving the other coordinates unchanged:
$$
(\tilde P_\mu f)(R_1,R_2,{\ldots} )=f(P_\mu(R_1,{\ldots} ,R_n),R_{n+1},{\ldots} ).
$$

\subsection{Examples of type I induced representations}

\subsubsection{Elementary representations}

Let us consider induced representations of type $\k=(\infty,k)$ induced
from Young subgroups with one infinite block and one finite block of size $n$.
First of all, note that all Young subgroups of this type are conjugate,
so that all representations of this type
are equivalent; thus we may fix an arbitrary partition $\Pi$ of type $(\infty,k)$,
say $\Pi=(\{1,{\ldots} ,k\},\{k+1,k+2,{\ldots}\})$, and speak of the
representation $I_{(\infty,k)}=I_\Pi$.
 By Theorem~\ref{th:typeI}(b), this representation is irreducible.

Recall the definition of an elementary representation of the infinite
symmetric group $\si$.

\begin{definition}
A representation of the infinite symmetric group $\si$
is called {\em elementary} if it is the inductive limit
of irreducible representations of finite symmetric groups $\sn$.
\end{definition}

\noindent{\bf Remark 1.}
Thus an elementary representation is realized in the $l^2$ space
on one class of tail-equivalent Young tableaux.

\smallskip
\noindent{\bf Remark 2.} In this definition, we do not require that
the embeddings preserve some sequence of fixed cyclic vectors.
However, it is easy to see that
the successive embeddings of the unit vector in the original
one-dimensional representation of $\S_1$ form such a sequence.

\medskip
Obviously, the representation $I_{(\infty,k)}$ of $\si$ is the inductive limit of the
representations $I_n$ of
the finite symmetric groups ${\mathfrak S}_n$ induced from the
identity representations
of the Young subgroups ${\frak S}_{\{1,{\ldots} ,k\}}\times
{\frak S}_{\{k+1,{\ldots},n\}}$. These representations are reducible; however,
it turns out that only one irreducible component ``survives'' in the limit,
so that the following proposition holds.

\begin{proposition}
\label{prop:elem}
The representation of the infinite symmetric group
induced from a Young subgroup of type $(k,\infty)$
is elementary.
\end{proposition}

\begin{proof}
Observe that a two-row diagram $\la^{(n)}=(n-k,k)$ is majorized in the dominance
ordering precisely by the diagrams $(n-m,m)$ with $m\le k$, and
all the corresponding Kostka numbers are equal to~$1$.
By~(\ref{ind}), we have $I_n=\sum_{m\le k}\pi_{(n-m,m)}$.
Consider the projection $P_{\la^{(n)}}\xi$
of the distinguished cyclic vector $\xi$
to the component $\pi_{\la^{(n)}}=\pi_{(n-k,k)}$. We have
$\dim I_n=\frac{n!}{k!(n-k)!}$, $\dim\la^{(n)}=\frac{n!(n-2k+1)}{k!(n-k+1)!}$,
and, obviously, $K_{\la^{(n)},\la^{(n)}}=1$.
Thus, by Lemma~\ref{l:spectral} (see Sec.~\ref{sec:spectral}),
$$
\|P_{\la^{(n)}}\xi\|^2=\frac{K_{\la^{(n)},\la^{(n)}}\dim\la^{(n)}}{\dim I_n}=
\frac{n-2k+1}{n-k+1}\to1\qquad\mbox{as }n\to\infty.
$$
It follows that the vectors $P_{\la^{(n)}}\xi$ converge to $\xi$, so that
the inductive limit of the irreducible
representations $\pi_{\la^{(n)}}$ coincides with $I_\Pi$, as required.
\end{proof}

It is well known (see, e.g., \cite{Nik, AAM, VTs} and also Sec.~\ref{sec:two})
that the representation $I_n$ can be realized in the space of symmetric tensors
of rank $k$ and dimension $n$. Then the representation $I_{(\infty,k)}$
can be realized in the space of infinite-dimensional symmetric tensors
of rank $k$.

\subsubsection{Representations of type $(\infty,\la)$}

Let us consider
induced representations of type $\k$ with
one infinite block and finitely many finite blocks, i.e.,
$k_\infty=1$ and
$k_1+k_2+{\ldots}=k<\infty$. Note that
all representations of this type are equivalent.

\begin{proposition}
Let $\Pi$ be a partition of $\N$ with one infinite block and finitely many finite blocks,
and denote by $\la=\la(\Pi)=1^{k_1}2^{k_2}{\ldots}$
the finite Young diagram of size  $k=k_1+k_2+{\ldots}$
formed by the lengths of the finite blocks of $\Pi$. Then
the decomposition if $I_\Pi$ into irreducible components is as follows:
\be
\label{elem}
I_\Pi=\sum_{|\mu|=k,\,\mu\unrhd\la}K_{\mu,\la}\pi_\mu^{1},
\ee
where  $\pi_\mu^{1}$ is
the elementary representation  of $\si$, namely, the inductive limit
$\pi_\mu^{1}=\lim_{n\to\infty}\pi_{\mu^{(n)}}$, where
$\mu^{(n)}=(n-k,\mu_1,\mu_2,{\ldots} )$
is the diagram with first row $n-k$ and the other rows forming the
diagram $\mu$.
\end{proposition}

Let us find the spectral measure of the distinguished cyclic
vector $\xi$ in this case. The arguments generalize the
proof of Proposition~\ref{prop:elem}.
Let $\la^{(n)}=(n-k,\la_1,\la_2,{\ldots} )$.
Fix a digram $\mu=(\mu_1,\mu_2,{\ldots})\unrhd\la$, $|\mu|=k$, set
$\mu^{(n)}=(n-k,\mu_1,\mu_2,{\ldots} )$, and let $n\to\infty$.
By Lemma~\ref{l:spectral},
$$
\|P_{\mu^{(n)}}\xi\|^2=\frac{K_{\mu^{(n)},\la^{(n)}}\dim\mu^{(n)}}{\dim I_n},
$$
where $I_n=\Ind_{\S_{\la^{(n)}}}^{\S_n}1$. (Note that
$I_\Pi$ is the inductive limit of $I_n$.)
We have $\dim I_n=\frac{n!}{(n-k)!\prod\la_i!}\sim \frac{n^k}{\prod\la_i!}$ as $n\to\infty$.
Further, it follows from the hook length formula that
$\dim\mu^{(n)}\sim\frac{n^k\dim\mu}{k!}$ as $n\to\infty$.
Finally, it can easily be seen from
the definition of Kostka numbers that $K_{\mu^{(n)},\la^{(n)}}=K_{\mu,\la}$
(indeed, $K_{\mu^{(n)},\la^{(n)}}$ is the number of fillings of
the shape $\mu^{(n)}$ with $n-k$ zeros, $\la_1$ ones, etc.; but, obviously,
the $n-k$ zeros must occupy exactly the first row of $\mu^{(n)}$).
Thus we obtain
\be
\label{cyl}
\lim_{n\to\infty}\|P_{\mu^{(n)}}\xi\|^2=\frac{K_{\mu,\la}\dim\mu\prod\la_i!}{k!}.
\ee
Note that the sum of the right-hand sides over all
$\mu\unrhd\la$ is equal to 1.
It follows that
the spectral measure of $\xi$ is discrete and supported by a finite set
of classes of tail-equivalent tableaux
indexed by diagrams $\mu\unrhd\la$. In particular,
$I_\Pi$ is a finite sum of elementary representations.

\medskip\noindent{\bf Remark 1.} Thus we see that all elementary representations
of $\si$ corresponding to inductive limits of sequences
of irreducible representations associated with Young diagrams
with growing first row can be
obtained by induction from Young subgroups.

\smallskip
\noindent{\bf Remark 2}.
Note that for $k_\infty\ge2$, the irreducible representation $\pi_\mu^{k_\infty}$ of $\si$
{\it is no longer elementary}.

\medskip{\bf Important particular case:
a hook with infinite hand and finite leg.} Let us consider the
representation
induced from ``a hook with infinite hand and finite leg,'' i.e.,
a Young subgroup with one infinite block and $n$ one-point blocks
(the case of a hook
with infinite hand and infinite leg will be considered
in Sec.~\ref{sec:infhook}, and that of a hook with
finite hand and infinite leg, in Sec.~\ref{sec:planch}).
In this case, $\la(\Pi)$ is the column diagram $1^n$ and
$\S_{\rm fin}=\{e\}$,
so that $I_{\rm fin}$ is the regular representation ${\rm Reg}_n$ of $\sn$. Thus
the decomposition (\ref{ind}) turns into the well-known formula
$$
{\rm Reg}_n=\sum_{|\mu|=n}\dim\mu\cdot\pi_\mu,
$$
and the decomposition (\ref{elem}) takes the form
\begin{equation}
\label{hookdec}
I_{(\infty,1^n)}=\sum_{|\mu|=n}\dim\mu\cdot\pi_\mu^{1}.
\end{equation}
In this case, the homogeneous space $X_\Pi$ can be identified
(by ``forgetting'' the infinite part of the partition $r$ in (\ref{hom}))
with the space
of $n$-sequences $(r_1,{\ldots} ,r_n)$ of distinct positive integers,
that is, $l^2(X_\Pi)$ can be identified with the space
of {\it diagonal-free infinite-dimensional tensors
of rank $n$}. On the other hand,
the regular representation $I_{\rm fin}={\rm Reg}_n$
has a realization in the space of tensors
of dimension $n$. The projections $P_\mu$
to the irreducible components of $I_{\rm fin}$
are the standard {\it Young symmetrizers}; and the projections to
the primary components, i.e.,
subspaces of tensors with given symmetry type, are
the {\it central Young symmetrizers} (see, e.g., \cite{Ham, Weyl}). Applying
these symmetrizers to
infinite-dimensional tensors of rank $n$, we obtain the irreducible (or primary)
components of $I_{(\infty,1^n)}$.

For example, the representation $I_{(\infty,1^2)}$ acts in the space
of infinite-dimensional tensors of rank 2, and
$$
I_{(\infty,1^2)}=\pi^1_{(2)}+\pi^1_{(1^2)},
$$
where the projections to the irreducible components are given by
symmetrizing and antisymmetrizing:
$$
P_{(2)}T_{r_1,r_2}=\frac12(T_{r_1,r_2}+T_{r_2,r_1}),\qquad
P_{(1^2)}T_{r_1,r_2}=\frac12(T_{r_1,r_2}-T_{r_2,r_1}).
$$

As follows from~(\ref{cyl}), the spectral measure of the distinguished cyclic vector
$\xi$ in this case has the Plancherel weights
(see Sec.~\ref{sec:markov}):
$$
\lim_{n\to\infty}\|P_{\mu^{(n)}}\xi\|^2=\frac{\dim^2\mu}{k!}.
$$

\subsubsection{Two-block representations}

By Theorem~\ref{th:typeI}(b), induced representations of type $\infty^2$,
induced from Young subgroups with two infinite blocks,
are irreducible. (Note that there is a continuum of nonequivalent
representations of this type.) As mentioned above, in the case of a two-row
diagram $\la$, all Kostka numbers $K_{\mu,\la}$ for $\mu\unrhd\la$
are equal to $1$. Thus we obtain
the following assertion.

\begin{proposition}
The representation of the infinite symmetric group
induced from a two-block Young subgroup
is irreducible, and its spectrum with respect to the
Gelfand--Tsetlin algebra is simple.
\end{proposition}

As shown in \cite{VTs}, in this case
one can obtain a complete spectral analysis of the induced
representations. We present these details in Sec.~\ref{sec:two}.

\section{Induced representations of type II}
\label{sec:factor}

\subsection{Small Young subgroups lead to representations
of type II}

Now let us consider small Young subgroups $\S_\Pi$.
We will assume that the partition $\Pi$ has finitely many
finite blocks of finite multiplicities and
denote by $\nu=\nu(\Pi)$
the Young diagram formed by the sizes of these blocks.

\begin{theorem}
\label{th:factor}
Let $\Pi$ be a partition of $\N$ of type
$\k=(\infty^{k_\infty},1^{k_1},2^{k_2},{\ldots} )$ that has finitely many
finite blocks of finite multiplicities
(i.e., $\sum_{j<\infty:\,k_j<\infty}k_j<\infty$), and assume that
there exists at least one $i\in\N$ with $k_i=\infty$. Then

{\rm(a)}
The representation $I_\Pi$ of the infinite symmetric group $\si$
is of type II.

{\rm (b)} This representation is a (type II) factor representation if and only if the
diagram $\nu(\Pi)$ consists of at most one row.

{\rm (c) (Central decomposition)}
Let $N=|\nu(\Pi)|$. The representation $I_\Pi$
is a finite sum of type II factor representations indexed by the primary
components of the representation $\Ind_{\S_{\nu(\Pi)}}^{\S_N}1$,
i.e., by Young diagrams $\mu$ such that $|\mu|=N$ and
$\mu\unrhd\nu(\Pi)$.
\end{theorem}

\begin{proof}
Let $J=J_\Pi=\{j\in\N:k_j=\infty\}$. By our assumptions, $J\ne\emptyset$.
For each $j\in J$, let $B^{(j)}_1,B_2^{(j)},{\ldots}$ be all blocks
of size $j$, and set $\S^{(j)}=\S_{B^{(j)}_1}
\times\S_{B_2^{(j)}}\times{\ldots} $.
In a similar way, let
$A_1,{\ldots} ,A_n$ be all finite blocks of finite multiplicities,
and let $\S_\nu=\S_{\nu(\Pi)}=\S_{A_1}
\times\S_{A_2}\times{\ldots}$ be the corresponding
finite Young subgroup of $\S_N$, where
$N=|A_1|+|A_2|+{\ldots} =|\nu(\Pi)|$. We also denote
$I_\nu=\Ind_{\S_\nu}^{{\frak S}_N} 1$.

It follows from the proof of Theorem~\ref{th:typeI} that the set of
intertwining operators for $I_\Pi$
(i.e., the commutant $\frak A'$ of the algebra $\frak A$ generated by the
representation operators)
decomposes into the tensor product
$$
\frak A'=R(I_\nu)\otimes\bigotimes_{j\in J}R_j,
$$
where $R(I_\nu)$ is the set of
intertwining operators for $I_\nu$ and $R_j$
is the algebra of operators generated by finite permutations of
the sets $B^{(j)}_1,B_2^{(j)},{\ldots}$. Obviously, $R_j$ is
(algebraically) isomorphic to the algebra generated by the regular
representation of the infinite symmetric group, which is
a factor of type~II$_1$. Thus $\otimes\bigotimes_{j\in J}R_j$
is a factor of type~II. Now if $\nu$ consists of a single row,
than the representation $I_\nu$ is irreducible and $R(I_\nu)$ consists
of scalar operators, so that the whole $\frak A'$
is a factor of type~II.
Otherwise, taking in $R(I_\nu)$ the projections to the primary
components of $I_\nu$, we obtain the central decomposition of $I_\Pi$
into a sum of factors.
\end{proof}

\noindent{\bf Remark.} Assume that $\#J=\#\{j\in\N:\,k_j=\infty\}<\infty$
and $\nu(\Pi)$ consists of at most one row, i.e., $I_\Pi$ is
a factor representation. Then the commutant $\frak A'=\otimes_{j\in J}R_j$
is a finite tensor product of factors of type~II$_1$, which is
a factor of type~II$_1$.
At the same time, it is not difficult to see that
if $\Pi$ is not the trivial partition into
singletons (i.e., $I_\Pi$ is not the regular representation),
then the algebra $\frak A$ itself
is a factor of type II$_\infty$.

\medskip
Theorems~\ref{th:typeI} and~\ref{th:factor}, which describe
induced representations of type~I and~II, respectively, do not
exhaust all induced representations. Namely, they leave out the case
when the partition $\Pi$ has  infinitely many finite blocks of finite multiplicities,
i.e., $\#\{i\in\N:\,k_i<\infty\}=\infty$. The most important example
of such a partition is a partition that has no infinite blocks
and at most one block of each finite size ($k_\infty=0$,
$k_i\le1$ for all $i\in\N$). {\it It is natural to conjecture  that the corresponding
representation is also of type~II.}

\subsection{Example: representations induced from a hook with
infinite hand and infinite leg}
\label{sec:infhook}

Let us consider in more detail the induced representation of
type $(\infty, 1^\infty)$ associated with a partition $\Pi$ having
one infinite block and infinitely many singletons.
Note that there is a continuum of nonequivalent representations of this type.

In this case, the homogeneous space $X$ is
the space of ``infinite-dimensional tensors
of infinite rank,'' i.e., infinite sequences
$(i_1,i_2,{\ldots} )$ of positive integers. It can be also
described as follows. Let
$\N_0$ be the subset of $\N$ consisting of all singletons of $\Pi$.
Then the space $l^2(X)$ can be identified with the space of
summable injections $f:\N_0\to\N$. The representation $I=I_\Pi$ is
generated by the left action of $\si$ by
$I_gf(n)=g^{-1}f(n)$. Denote by $\frak A$ the von Neumann
algebra generated by this representation. It is easy to see that the
commutant ${\frak A}'$ of this algebra, i.e., the set of intertwining
operators, is generated by the right action of the group
$\S(\N_0)$ of finite permutations of $\N_0$ by substitutions:
$T_{\sigma} f(n)=f(\sigma^{-1}n)$. In particular, ${\frak A}'$ is algebraically
isomorphic to the von Neumann algebra generated by the regular
representation of the infinite symmetric group, i.e., is a factor of type
II$_1$. As to the factor $\frak A$ itself, it is easy to see that
it decomposes into the infinite direct sum of factors of type II$_1$.
Indeed,
for each infinite subset $B\subset\N$,
let $H_B$ be the subspace of $l^2(X)$
consisting of functions $f$ such that $f^{-1}(\N_0)=B$.
This subspace is obviously invariant with respect to ${\frak A}'$, so that
$l^2(X)=\oplus_{B}H_B$
is the desired decomposition of $\frak A$
into the infinite direct sum of factors of type II$_1$.
Thus we obtain the following proposition.

\begin{proposition}
The induced representation of
type $(\infty, 1^\infty)$ of the infinite symmetric group $\si$
is a factor representation. The von Neumann algebra $\frak A$  generated by
the representation operators is a factor of type II$_\infty$,
a direct sum of isomorphic factors of type II$_1$.
Its commutant ${\frak A}'$ is a factor of type II$_1$ algebraically
isomorphic to the factor generated by the regular representation of $\si$.
\end{proposition}

Note that the distinguished cyclic
vector $\xi$, which in this realization is
just the identical function $\xi(n)\equiv n$, is a cyclic
vector for $\frak A$, but not for ${\frak A}'$.
Thus it defines a finite trace on the commutant
by the formula $\tr(A)=(A\xi,\xi)$, $A\in{\frak A}'$
(which is just the trace corresponding to the regular representation
of the group $\S(\N_0)$, i.e., coincides on the elements of $\S(\N_0)$
with the delta function at the identity element),
but there is no finite trace on $\frak A$.

\medskip\noindent{\bf Remark.}
The described construction is a particular case of the following one.
We have the standard left action of a group $G$ on the discrete homogeneous space
$X=G/H$, where $H$ is a subgroup of $G$, and the corresponding
unitary representation of $G$  in
$l^2(X)$. The commutant of the corresponding algebra is generated
by the automorphisms of the $G$-space $X$, i.e., by
the right action of the group
$N(H)/H$, where $N(H)$ is the normalizer
of $H$ in $G$. The specific property of this situation is
that both left and right representations are generated by
substitutional actions  of $G$ and $N(H)/H$, respectively (cf.~\cite{strange}).

\section{Examples of the spectral analysis of induced representations}
\label{sec:spectral}

In this section, we will consider two nontrivial examples
of classes of induced representations for which we can explicitly find
the spectral measure of the distinguished cyclic vector and,
consequently, obtain a spectral realization.

To this end, the following simple lemma is useful.
Consider the restriction $I_n$ of $I_\Pi$ to a finite symmetric group $\sn$
in the cyclic hull $\sn\xi$ of the distinguished cyclic vector $\xi$
(which is the representation induced from the identity representation of
the Young subgroup $H=\S_{\Pi_n}$ of $\sn$ associated with the partition
$\Pi_n=\Pi\cap\{1,{\ldots} ,n\}$). Let $\la$ be a Young diagram
with $n$ cells, and denote by $P_\la$ the projection from $I_n$
to the {\it primary} component $V_\la$
that is a multiple of $\pi_\la$.

\begin{lemma}
\label{l:spectral}
The squared norm
$\|P_\la\xi\|^2$ is equal to the relative dimension of $V_\la$
in $I_n$, i.e.,
\be
\label{primproj}
\|P_\la\xi\|^2=\frac{\dim V_\la}{\dim I_n}.
\ee
\end{lemma}
\begin{proof}
It is well known (see, e.g., \cite{Serre})
that the projection $P_\la$ is given by the formula
$$
P_\la=\frac{\dim\la}{n!}\sum_{g\in\sn}\chi_\la(g)I_n(g),
$$
where $\chi_\la$ is the character of $\pi_\la$. Then
$$
(P_\la\xi,P_\la\xi)=(P_\la\xi,\xi)=\frac{\dim\la}{n!}\sum_{g\in\sn}
\chi_\la(g)(I_n(g)\xi,\xi)=\frac{\dim\la}{n!}\sum_{h\in H}\chi_\la(h),
$$
because $(I_n(g)\xi,\xi)$ is equal to $1$ if $g\in H$ and $0$ otherwise.
But the latter sum is equal to
$(\chi_\la|_H,1_H)\cdot|H|$,
which is in turn equal to $(\chi_\la, I_n)\cdot|H|$
by  the Frobenius reciprocity, so that we obtain
$$
\|P_\la\xi\|^2=\frac{(\chi_\la, I_n)\dim\la\cdot|H|}{n!}=
\frac{\dim V_\la}{\dim I_n},
$$
because $(\chi_\la, I_n)$ is the multiplicity of $\pi_\la$ in $I_n$
and $\dim I_n=\frac{n!}{|H|}$.
\end{proof}

\subsection{Markov representations of the infinite symmetric group}
\label{sec:markov}

In this section, we  recall necessary notions from the representation
theory of the symmetric groups and the notions of simple and Markov
representations.

Denote by $T_\la$ the set (consisting of $\dim\la$ elements)
of Young tableaux of shape
$\la\in\yn$, or, which is the same, the set of paths in the Young graph
from the empty diagram $\emptyset$ to $\la$.
Let
$T_n=\cup_{\la\in\yn}T_\la$
be the set of Young tableaux with $n$ cells.
According to the branching rule for irreducible representations
of the symmetric groups, the space $V_\la$ of the irreducible
representation  $\pi_\la$ decomposes into the sum of one-dimensional subspaces
indexed by the tableaux $u\in T_\la$. The basis $\{h_u\}_{u\in T_\la}$
consisting of vectors of these subspaces is called the
{\it Gelfand--Tsetlin basis}. It is an eigenbasis for the {\it Gelfand--Tsetlin
algebra} $\gzn$, the subalgebra in the group algebra $\gn$ generated by the centers
$Z[{\mathfrak S}_1], Z[{\mathfrak S}_2],{\ldots} ,Z[{\mathfrak S}_n]$
(see \cite{VO}).

Denote by $T=\varprojlim T_n$ the space of infinite Young tableaux
(the projective limit of $T_n$ with respect to the natural projections
forgetting the tail of a path). With the topology of coordinatewise convergence
$T$ is a totally disconnected metrizable compact space.
The tail equivalence relation $\sim$ on $T$ is defined as follows:
paths $s=(s_1,s_2,{\ldots})$ and
$t=(t_1,t_2,{\ldots})$ are equivalent
if and only if
$s_k=t_k$ for all sufficiently large $k$.
Denote by $[t]_n\in T_n$ the initial segment of length $n$
of a tableau $t\in T$. Given a finite tableau $u\in T_n$, denote
by $C_u=\{t:[t]_n=u\}$ the corresponding cylinder set;
for $\la\in\yn$, let $C_\la=\{t:t_n=\la\}=\cup_{u\in T_\la}C_u$
be the set of all paths passing through $\la$.

\begin{definition}
A measure $M$ on the space $T$ is called {\em Markov}
if for every $n\in\N$ the following condition holds:
for any diagrams $\la\in\yn$ and
$\La\in{\mathbb Y}_{n+1}$ such that $\La\subset\la$
and for any path $u\in T_\la$, the events
$C_u$ (``the past'') and $C_\La$ (``the future'') are independent
given
$C_\la$ (``the present'').
\end{definition}

In other words,
a random tableau $t=(t_1,t_2,{\ldots} )$,
regarded as a sequence of random variables $t_n$, where $t_n$
takes values in the set $\yn$ of Young diagrams with $n$ cells,
is a Markov chain in the ordinary sense. In terms of
transition probabilities, this means that
the transition probability $\frac{M(C_\La\cap C_u)}{M(C_u)}$
depends only on the form $\la$ of a tableau $u$, but not on the tableau itself.
Note that the ``forward'' and ``backward'' Markov properties are
equivalent, so that the definition of a Markov measure can be formulated
in a similar way in terms of cotransition probabilities.

One of the most important examples of Markov measures on $T$ is the
{\it Plancherel measure} $P$, which is the spectral measure of the
regular representation of $\si$.

\begin{definition}
The Plancherel measure on the space $T$ of infinite Young tableaux
is the Markov measure with transition probabilities
$\Prob(\la,\La)=\frac{\dim\La}{(n+1)\dim\la}$, where $\la\in\mathbb Y_n$,
$\La\in\mathbb Y_{n+1}$, $\la\subset\La$.
\end{definition}

The cylinder distributions $P_n$ of the Plancherel measure
are given by the formula $P_n(C_\la)=\frac{\dim^2\la}{n!}$, $\la\in\mathbb Y_n$.
This distribution on the set of Young diagrams with $n$ cells
is called the {\it Plancherel measure} on Young diagrams
(see \cite{77}).

Now assume that we have a cyclic representation $\pi$ of the group
$\sn$ in a space $V$ that has a simple spectrum (i.e.,
decomposes into the sum of pairwise nonequivalent irreducible representations)
and a unit cyclic vector $\xi$ in this representation.

\begin{definition}
We say that $\xi$ is a {\it Markov} vector
if its spectral measure with respect to the Gelfand--Tsetlin algebra
is Markov.
\end{definition}

\begin{lemma}
\label{l:mark}
Let $\pi$ be a unitary representation of the group
$\sn$ with simple spectrum. A cyclic vector
$\xi$ of the representation $\pi$ is Markov if and only if
for every $k<n$, the representation of the group
${\mathfrak S}_k$ in the cyclic hull
${\mathfrak S}_k\xi$ of $\xi$ with respect to
${\mathfrak S}_k$ has a simple spectrum.
\end{lemma}

\begin{proof}
Clearly, a cyclic vector is Markov if for every
$k<n$ and every diagram $\la\in{\mathbb Y}_k$, the probability
of any tableau with this diagram does not depend
on the continuation of this tableau to the level $n$.
In terms of representations and cyclic vectors, this means
that the norm of the projection of the cyclic vector to
the subspace of the representation of the group
${\mathfrak S}_k$ equivalent to $\pi_\la$
does not depend on the way in which we have arrived at this subspace.

Now we use the following simple lemma.

\begin{lemma}
\label{l:cyclic}
Assume that in a finite-dimensional Hilbert space
$H$ there is a unitary representation of a group $G$ that
is primary, i.e., decomposes into the direct
(not necessarily orthogonal) sum
$H=H_1\oplus H_2\oplus{\ldots} \oplus H_n$
of equivalent irreducible representations, and in each of them
there is a cyclic vector
$v_i\in H_i$, $i=1,{\ldots} ,n$.
Then the following two assertions are equivalent:

{\rm 1.} For any $i,j$, there exists an isometry
$T_{i,j}:H_i\to H_j$ intertwining the corresponding representations
such that $T_{i,j}v_i=v_j$.

{\rm 2.} In the cyclic hull of the vector $v=\sum v_i$,
the representation is irreducible.
\end{lemma}

\begin{proof}
Without loss of generality assume that $n=2$. Then the representation
in the cyclic hull
of $v$ is irreducible if and only if this vector is of rank $1$, i.e.,
$v$, regarded as an element of the tensor product, has the form
$x\otimes x$.
\end{proof}

Now let $\xi_\mu$ be the projection of $\xi$ to the irreducible
component $\pi_\mu$ of $\pi$. Fix $k<n$ and an irreducible
representation $\pi_\la$ of $\S_k$. Consider
the restriction of $\pi$ to $\S_k$, and let
$H=H_1\oplus H_2\oplus{\ldots} \oplus H_n$, where $H_i\simeq\pi_\la$,
be the primary component of this restriction corresponding to $\pi_\la$.
The Markov property means that the norms of the projections
of $\xi_\mu$ to $H_i$ coincide, which is equivalent, in
view of Lemma~\ref{l:cyclic}, to the fact that $\pi_\la$ has multiplicity $1$
in the decomposition
of the representation of ${\mathfrak S}_k$ in the cyclic hull
${\mathfrak S}_k\xi$ of $\xi$.
\end{proof}

Now let us consider representations of the infinite symmetric group $\si$.

If we are given a quasi-invariant measure $\mu$ on the space of Young tableaux
$T$ and a $1$-cocycle $c$ on pairs of tail-equivalent paths
taking values in the group of complex numbers of modulus $1$,
then we can construct a unitary representation of the group $\si$
in the space $L^2(T,\mu)$ as follows (see, e.g., \cite{Itogi}). Recall
that the Fourier transform allows one to realize the group algebra
$\C[\si]$ of the infinite symmetric group as the cross product
constructed from the commutative algebra of functions on the space of tableaux $T$
(Gelfand--Tsetlin algebra) and the tail equivalence relation.
The desired representation is given by
\be
L_g h(s)=\sum_{t\sim s}\sqrt{\frac{d\mu(s)}{d\mu(t)}}\hat g(s,t)c(s,t)h(t),
\qquad h\in L^2(T,\mu),
\label{repr}
\ee
where $\hat g$ is the function on pairs of tail-equivalent paths
corresponding to an element $g\in\si$
(the Fourier transform of $g$). Note that the cocycle is trivial on
the space of finite tableaux.

\begin{definition}
A representation of the infinite symmetric group
$\si$ is called
{\em simple} if it is the inductive
limit of representations of the finite symmetric groups
$\sn$ with simple spectrum.
\end{definition}

\begin{definition}
A representation $\pi$ of the group
$\si$ with simple spectrum  is called
{\em Markov} if the space of $\pi$ contains a cyclic vector
whose spectral measure (with respect to the Gelfand--Tsetlin algebra)
is Markov.
\end{definition}

Note that a representation with simple spectrum is Markov if and only if
the measure $\mu$ in its realization~(\ref{repr}) is Markov.

The following theorem is an easy consequence of Lemma~\ref{l:mark}.

\begin{theorem}[\cite{VTs}]
\label{th:mark}
A representation of the infinite symmetric group
is Markov if and only if it is simple.
\end{theorem}

\subsection{The spectral analysis of two-block induced representations}
\label{sec:two}

In this section, we present the complete spectral analysis
of induced representations
of type $\infty^2$.
By Theorem~\ref{th:typeI}, these
representations are irreducible.
It turns out that in this case the spectral measure
of the distinguished cyclic vector is a Markov measure on the space
of infinite Young tableaux $T$.
The material of this section is mostly borrowed from the authors' paper
\cite{VTs}. Theorem~\ref{th:spectral} is proved there using
the so-called {\it tensor model} of two-row representations
of the symmetric groups (i.e., representations induced
from two-block Young subgroups).
See~\cite{VTs} for a detailed
description of this model (in particular, explicit formulas
for the Gelfand--Tsetlin basis).

A partition $\mathbb N=A\cup B$ of type $\infty^2$
is uniquely determined  by
an infinite sequence
$\xi=\xi_1\xi_2{\ldots}$ of $0$'s and
$1$'s (an ``infinite tensor''),
where $\xi_i=1$ if $i\in A$, and $\xi_i=0$ if $i\in B$.
Then the induced representation in question is equivalent to
the natural substitutional representation of $\si$
on infinite sequences in the cyclic hull
of the sequence $\xi$, which we will denote
by $\pi_\xi$.
Note that the orbit of $\xi$ is the discrete set $O_\xi$ of infinite sequences
of $0$'s and $1$'s eventually coinciding with $\xi$,
and $\pi_\xi$ is a unitary representation of $\si$
in the space $l^2(O_\xi)$.

For simplicity, it is convenient to assume that the number of
$1$'s among the first $n$ elements of $\xi$ does not exceed $n/2$.
It is not difficult to see that an arbitrary case can be reduced
to this one, but we omit the corresponding technical details.

Denote by $\la_{n,k}=(n-k,k)$ the diagram with two rows of lengths
$n-k$ and $k$.

\begin{theorem}
\label{th:spectral}
The spectral measure $\mu_\xi$
of the cyclic vector $\xi$ in the representation $\pi_\xi$
with respect to the Gelfand--Tsetlin algebra is a Markov measure on the
space of infinite Young tableaux $T$, and its transition probabilities
are given by the
following formula. Denote by $m(n)$ the number of $1$'s among the first
$n$ elements of $\xi$.

If $\xi_{n+1}=0$, then
\be
\label{8}
\Prob(\la_{n,k},\la_{n+1,k})=\frac{n-m(n)-k+1}{n-2k+1},\qquad
\Prob(\la_{n,k},\la_{n+1,k+1})=\frac{m(n)-k}{n-2k+1}.
\ee

If $\xi_{n+1}=1$, then
\be
\label{9}
\Prob(\la_{n,k},\la_{n+1,k})=\frac{m(n)-k+1}{n-2k+1},\qquad
\Prob(\la_{n,k},\la_{n+1,k+1})=\frac{n-m(n)-k}{n-2k+1}.
\ee
\end{theorem}

Note that all spectral measures of induced representations
considered in Theorem~\ref{th:spectral}
are not central (except for the trivial case when one of the sets
in the partition is empty).

\smallskip\noindent{\bf Example.}
Let $\xi=0101{\ldots}$. Then $m(n)=[n/2]$, and the formulas for transition
probabilities take the following form:

$\bullet$ if $n$ is odd,
\be
\Prob(\la_{n,k},\la_{n+1,k})=\frac{n-2k+2}{2(n-2k+1)},\qquad
\Prob(\la_{n,k},\la_{n+1,k+1})=\frac{n-2k}{2(n-2k+1)};
\label{12odd}
\ee

$\bullet$ if $n$ is even,
\be
\Prob(\la_{n,k},\la_{n+1,k})=
\Prob(\la_{n,k},\la_{n+1,k+1})=\frac12.
\label{12even}
\ee

It is convenient to rewrite formulas~(\ref{12odd}), (\ref{12even})
introducing the change of indices $j=n-2k$. In these terms,
a Young tableau is determined by a sequence
$(j_1,j_2,{\ldots})$, where  $j_n$ takes the values
$0,1,{\ldots} ,n$,
and the transition probabilities of the measure $\mu_\xi$ are equal to
$$
\Prob(j,j+1)=\frac{j+2}{2(j+1)},\qquad
\Prob(j,j-1)=\frac{j}{2(j+1)}
$$
at an odd moment of time; and
$$
\Prob(j,j+1)=\Prob(j,j-1)=\frac12
$$
at an even moment of time.
We see that a random Young tableau governed by the measure $\mu_\xi$
is a trajectory of
a nonhomogeneous (neither in time nor in space)
random walk on ${\mathbb Z}_+$.
Thus the induced representations of the infinite symmetric group
considered in this paper act in spaces of functions over
trajectories of natural random walks. Explicit formulas for
this action are given by Young's orthogonal form.

\subsection{Spectral measure of representations of type $(1^\infty,\nu)$}
\label{sec:planch}

Consider the representations $I_\Pi$ induced from partitions $\Pi$ of type  $(1^\infty,\nu)$
with no infinite blocks, infinitely many singletons, and finitely
many finite blocks of length greater than one. Here, as above, we denote
by $\nu$ the finite Young diagram formed by the lengths of finite blocks
of finite multiplicities, and let $|\nu|=n$.

According to Theorem~\ref{th:factor}, $I_\Pi$ decomposes into a finite
sum of factor representations $\rho_\mu$ indexed by Young diagrams
$\mu$ with $n$ cells such that $\mu\unrhd\nu$.
Let us find the spectral measure $M=M^\xi$ of
the distinguished cyclic vector $\xi$ with respect to the
Gelfand--Tsetlin algebra.

Given $N\ge n$,
denote by $\nu_N$ the diagram obtained from $\nu$ by adding
$N-n$ rows of length $1$. Let $\la$ be a Young diagram with $N$ cells.
By formula~(\ref{ind}), the definition of the spectral measure,
and Lemma~\ref{l:spectral}, the cylinder distribution $M_N$ of $M$
is given by
$$
M_N(C_{\la})=\|P_{\la}\xi\|^2
=\frac{\prod\nu_i!}{N!}K_{\la,\nu_N}\dim\la.
$$
It is not difficult to see from the definition of the  Kostka numbers that
$$
K_{\la,\nu_N}=\sum_{\mu\unrhd\nu}\dim(\mu,\la)K_{\nu,\mu},
$$
where $\dim(\mu,\la)$ is the number of paths in the Young graph
from $\mu$ to $\la$. Therefore, we obtain
$$
M_N(C_{\la})=\sum_{\mu\unrhd\nu}\frac{\dim(\mu,\la)K_{\nu,\mu}\dim\la\prod\nu_i!}{N!}=
\sum_{\mu\unrhd\nu}\frac{K_{\nu,\mu}\dim\mu\prod\nu_i!}{n!}\cdot
\frac{n!\dim(\mu,\la)\dim\la}{\dim\mu\cdot N!}.
$$
But it is not difficult to check that the second quotient in the right-hand side
of the latter formula is exactly the conditional distribution
$P(\cdot\,|\,t_n=\mu)$ of the
Plancherel measure $P$ on the space $T$ of infinite Young tableaux
$t=(t_1,t_2,{\ldots} )$
(paths in the Young graph) given that at the $n$th level $t$ passes
through $\mu$. Note also that the first quotient is the relative dimension
of the primary representation corresponding to the diagram $\mu$ in
$\Ind_{\S_\nu}^{\sn}$, i.e., $M_{\rm fin}(\mu)$, where $M_{\rm fin}$ is
the spectral measure of the distinguished cyclic vector in the
representation $\Ind_{\S_\nu}^{\sn}$ of the finite symmetric group $\sn$.
Thus we obtain the following result.

\begin{proposition}
\label{prop:sp}
The spectral measure $M=M^\xi$ of the distinguished cyclic vector $\xi$
in the representation $I_\Pi$
is a convex combination of
conditional Plancherel measures:
$$
M=\sum_{\mu\unrhd\nu}M_{\rm fin}(\mu)\cdot P(\cdot\,|\,t_n=\mu),
$$
where $M_{\rm fin}(\mu)=\frac{K_{\nu,\mu}\dim\mu\prod\nu_i!}{n!}$
is the spectral measure of the distinguished cyclic vector in the
representation $\Ind_{\S_\nu}^{\sn}$ of the finite symmetric group $\sn$.
In particular, it is absolutely continuous with respect to the Plancherel
measure $P$ with piecewise constant (cylinder) density
$$
\frac{dM}{dP}(t)=\frac{K_{\nu,\mu(t)}\prod\nu_i!}{\dim\mu(t)}\quad\mbox{if }
t\in T\mbox{ passes through the diagram } \mu(t)\mbox{ at level } n.
$$
\end{proposition}

The Plancherel measure is a central Markov measure
on the space of infinite Young diagrams. However,
the spectral measure $M$ is not central and is not Markov. However, it is
{\it multi-Markov}, in the sense that if we ``glue'' the first $n$
levels in one block, i.e., consider a random Young tableau as a
sequence of Young diagrams $(\square, \la_n,\la_{n+1},{\ldots} )$, then
it will be a Markov chain.

The following example is the simplest case of a representation of type
$(1^\infty,\nu)$, in which the representation is a factor.

\medskip\noindent{\bf Example. A hook with finite hand and infinite leg.}
If $\nu=(n)$, i.e., the partition consists of one finite block
of size $n$ and infinitely many singletons (``a
hook with finite hand and infinite leg''), then
by Theorem~\ref{th:factor}(b) the induced representation
is a factor, and by Proposition~\ref{prop:sp}
the spectral measure of the distinguished cyclic
vector is the conditional distribution $P(\cdot\,|\,t_n=(n))$ of the Plancherel measure
given that at the $n$th level $t$ passes through the
one-row diagram $(n)$.


\begin{thebibliography}{55}

\bibitem{Binder1}
M.~W.~Binder, Irreducible induced representations of ICC-groups,
{\it Math. Ann.} {\bf 294}  (1992),  37--47.

\bibitem{Binder2}
M.~W.~Binder, On induced representations of discrete groups,
{\it Proc. Amer. Math. Soc.}  {\bf118}  (1993),  301--309.

\bibitem{Binder3}
M.~W.~Binder, Induced factor representations of discrete groups and
their types, {\it J. Funct. Anal.} {\bf115} (1993), 294--312.

\bibitem{Fulton}W.~Fulton, {\it Young Tableaux. With Applications to
Representation Theory and Geometry}. Cambridge University Press, Cambridge, 1997.

\bibitem{Ham} M.~Hamermesh, {\it Group Theory and its Application to Physical
Problems}, Dover Publ., New York, 1989.

\bibitem{Hirai1} T.~Hirai, Some aspects in the theory of representations
of discrete groups. I, {\it  Proc. Japan Acad. Ser. A Math. Sci.}
{\bf66}  (1990),  315--318.

\bibitem{Hirai2} T.~Hirai,
Construction of irreducible unitary representations of the infinite
symmetric group $S_\infty$, {\it J. Math. Kyoto Univ.}
{\bf31}  (1991),  495--541.

\bibitem{Mac}
I.~Macdonald, {\it Symmetric Functions and Hall Polynomials}, 2nd edition.
Clarendon Press, Oxford, 1995.

\bibitem{Mackey} G.~Mackey, {\it The Theory of Unitary Group Representations},
The University of Chicago Press, Chicago--London, 1976.

\bibitem{Nik} P.~P.~Nikitin, A realization of the
irreducible representations of $S_n$
corresponding to 2-row diagrams in
square-free symmetric multilinear forms. {\it Zapiski Nauchn. Semin. POMI}
{\bf 301} (2003), 212--219. English translation:
{\it J. Math. Sci. (New York)}  {\bf129}, No.~2, 3796--3799 (2005).

\bibitem{Obata1} N.~Obata, Certain unitary representations of the
infinite symmetric group, I, {\it Nagoya Math. J.} {\bf 105} (1987), 109--119.

\bibitem{Obata2} N.~Obata,
Some remarks on induced representations of infinite discrete groups,
{\it Math. Ann.}  {\bf 284}  (1989), 91--102.

\bibitem{Serre} J.-P.~Serre, {\it Linear Representations of Finite
Groups}, Springer-Verlag, New York--Heidelberg, 1977.

\bibitem{AAM}
N.~V.~Tsilevich and A.~M.~Vershik, On different models
of representations of the infinite symmetric group.
To appear in {\it Adv. Appl. Math.}


\bibitem{strange} A.~M.~Vershik, Strange factor representations of type
II$_1$, and pairs of dual dynamical systems.
{\it Moscow Math. J.} {\bf 3}, No. 4 (2003), 1141--1157.

\bibitem{V} A.~M.~Vershik,
Inductive proof of Young's rule
(New approach to the representation theory of the symmetric groups. III),
to appear in {\it Moscow Math. J.}

\bibitem{77} A.~M.~Vershik and  S.~V.~Kerov,
Asymptotics of the Plancherel measure of the symmetric group and the limiting form of Young tableaux,
{\it Dokl. Akad. Nauk SSSR} {\bf233}, No.~6, 1024--1027 (1977).
English translation: {\it Sov. Math. Dokl.} {\bf18}, 527--531 (1977).

\bibitem{Itogi}  A.~M.~Vershik and  S.~V.~Kerov,
Locally semisimple algebras.  Combinatorial theory and
   the $K$-functor. In:
   {\it Itogi Nauki i Tekhniki, Ser. Sovrem. Probl. Mat.}, Vol. 26.
   VINITI, Moscow, 1985, pp.~3--56.
   English translation: {\it J. Sov. Math.} {\bf38}, 1701--1733 (1987).

\bibitem{GB} A.~M.~Vershik and S.~V.~Kerov, The Grothendieck group of
infinite symmetric group and symmetric functions (with the elements of
the theory of K$_0$-functor of AF-algebras). In:
A.~M.~Vershik and D.~P.~Zhelobenko (eds.), {\it Representation of Lie
Groups and Related Topics},
Gordon and Breach Sci. Publ., 1990, pp. 39--118.

\bibitem{VO} A.~M.~Vershik and A.~Yu.~Okounkov,
A new approach to the representation theory of symmetric groups. II.
{\it Zapiski Nauchn. Semin. POMI} {\bf 307}, 57--98 (2004).
English translation: {\it J. Math. Sci. (New York)} {\bf 131}, No.~2,
5471--5494 (2005).

\bibitem{325} A.~M.~Vershik and N.~V.~Tsilevich,
On the Fourier transform on the infinite symmetric group.
{\it Zap. Nauchn. Semin. POMI} {\bf 325} (2005), 61--82.
English translation to appear in {\it J. Math. Sci. (New York)}.

\bibitem{VTs}A.~M.~Vershik and N.~V.~Tsilevich,
Markov measures on Young tableaux
and induced representations of the infinite symmetric group,
{\it Probab. Theory Appl.} {\bf51} (2006), 47--63.

\bibitem{Weyl}H.~Weyl, {\it The Theory of Groups and Quantum Mechanics},
Dover Pub., New York, 1949.


\end{thebibliography}
\end{document}